\documentclass[12pt]{amsart}
\usepackage{graphicx}
\usepackage{amsrefs}
\usepackage{amssymb}
\usepackage[all]{xy}
\SelectTips {cm}{}
\usepackage{pinlabel}
\usepackage{pdfsync}
\usepackage{amsmath}
\usepackage{amsfonts}
\usepackage{amsthm}
\usepackage{etoolbox}
\usepackage{hyperref}

\makeatletter
\def\@MMRR#1 #2@{\hbox{\ \href{https://mathscinet.ams.org/mathscinet-getitem?mr=#1}{MR}\hskip2pt\@fullMR}}
\def\MR#1{\def\@fullMR{#1}\expandafter\@MMRR#1 M@\def\xx{#1}}

\newcounter{localNumber}
\newcommand{\resetSectionLabel}{\setcounter{localNumber}{0}}%

\newcommand{\nextLabel}[2][]{\addtocounter{localNumber}{1}\def\LRTtwo{#2}\ifstrempty{#1}%
{\let\@nextLabel\@@nextLabelEmpty}%
{\def\LRTone{#1}\let\@nextLabel\@@nextLabel}\@nextLabel}

\def\@@nextLabelEmpty{\expandafter\ifdefempty{\LRTtwo}{}{\LRTtwo\ }\number\value{section}.\number\value{localNumber}}

\def\@@nextLabel{%
\expandafter\xdef\csname\LRTone\endcsname{{\number\value{section}.\number\value{localNumber}}}%
\immediate\write\@auxout{\string\expandafter\string\gdef\string\csname\space\LRTone\string\endcsname{\csname\LRTone\endcsname}}%
\expandafter\xdef\csname\LRTone Link\endcsname{\noexpand\hyperlink{\LRTone}{\csname\LRTone\endcsname}}%
\immediate\write\@auxout{\string\expandafter\string\gdef\string\csname\space\LRTone Link\string\endcsname{\noexpand\hyperlink{\LRTone}{\csname\LRTone\endcsname}}}%
\expandafter\xdef\csname\LRTone Name\endcsname{\LRTtwo}%
\immediate\write\@auxout{\string\expandafter\string\gdef\string\csname\space\LRTone Name\string\endcsname{\LRTtwo}}%
\expandafter\ifdefempty{\LRTtwo}{}{\LRTtwo\ }%
\hypertarget{\LRTone}{\number\value{section}.\number\value{localNumber}}%
}

\theoremstyle{plain}

\newtheorem*{Theorem*}{\@thmName}

\theoremstyle{definition}

\newtheorem*{Definition*}{\@defName}

\gdef\empty@LRTU{}

\newenvironment{DefS}[2][]{\def\@defName{#2}%
\begin{Definition*}}%
{\end{Definition*}}

\newenvironment{ThmS}[2][]{\def\@thmName{#2}%
\begin{Theorem*}}%
{\end{Theorem*}}
\makeatother
\def\xxAfLink{}
\def\oDqLink{}
\def\oDrLink{}

\def\smfrac#1#2#3{{\frac{\scriptscriptstyle#1}{\raise#3pt\hbox{$\scriptscriptstyle#2$}}}}
\def\conj#1{\overline{#1}}
\def\MOD#1{\hbox{\rm mod $#1$}}
\def\rank#1{{\rm rank}\>#1}
\def\preT#1#2{{_{#1}{#2}}}
\def\dualT#1#2{{#1}\cdot{#2}}
\def\ad{Ad}
\def\qs#1{#1^{s}}
\def\ir#1#2{#1_{#2^\perp/#2}}
\def\Pont{Pontrjagin}
\def\dual#1{\textrm{Hom}(#1,\qmodz)}
\def\qh#1{#1^q}
\def\b{b}
\def\R{\mathbb R}
\def\Z{\mathbb Z}
\def\Q{\mathbb Q}
\def\qmodz{\Q/\Z}
\def\rmodz{\R/\Z}
\def\cy#1{\Z/#1\Z}
\def\q{\psi}
\def\c{\omega}
\def\arf{\beta}
\def\sig{\sigma}
\def\cupp{\>{\scriptstyle\cup}\> }
\def\rof#1{\langle#1\rangle}
\def\p{{\mathcal P}}
\def\abf#1{b_{#1}}
\def\exp#1{e^{2\pi i{#1}}}
\def\magg#1{N(#1)}
\def\wq#1{{\mathcal W}(#1)}
\def\enh#1{\Gamma(#1)}
\def\red#1{{#1}^{red}}
\def\m{\hbox to0pt{\hss$-$}}
\def\nqr{{\eta}}
\def\gnqr{{\gal{\nqr}}}
\def\qrN#1(#2){\ell_{#1}(#2)}
\def\gal#1{\gamma_{#1}}
\def\uf#1{\textbf{1}_{#1}}
\def\IM{{\rm Im}}

\begin{document}
\title{Gauss Sums in Algebra and Topology}
\author{Laurence R. Taylor}
\email{taylor.2@nd.edu}
\thanks{Partially supported by the N.S.F.}
\address{
Department of Mathematics
University of Notre Dame
Notre Dame, IN 46556 
}

\begin{abstract}
We consider Gauss sums associated to functions $T\to\rmodz$
which satisfy some sort of ``quadratic'' property and investigate
their elementary properties.
These properties and a Gauss sum formula from the
nineteenth century due to Dirichlet give 
the Milgram Gauss sum formula computing the signature 
mod $8$ of a non{--}singular bilinear form over $\Q$.
Brown derived some results on the signature mod $8$ 
of non{--}singular integral forms.
Kirby and Melvin gave a formula for a generalization of this
invariant to possibly non{--}singular forms and we further
generalize it here.
The Milgram Gauss sum formula and these formulas allow 
us to reprove Brown's result
without resort to Witt group calculations.
Assuming a bit of algebraic topology, we  reprove
a theorem of Morita's computing the signature
mod $8$ of an oriented Poincar\'e duality space from the 
\Pont\ square without using Bockstein spectral sequences.
Since we work with forms which may be singular, we also 
obtain a version of Morita's theorem for Poincar\'e spaces 
with boundary.
Finally we apply our results to the bilinear form
$Sq^1x\cup y$ on $H^1(M;\cy2)$ of an orientable $3${--}manifold.
\end{abstract}
\maketitle

\section{Introduction.}\resetSectionLabel

Gauss sums have a long and venerable history. 
A general version has a finite set $T$, a function
$\q\colon T\to\rmodz$ and the associated 
Gauss sum
\[G(\q)=\sum_{t\in T} \exp{\q(t)}\ .\]
Even more general notions would replace $T$ by
a measure space and the finite sum by an integral.
The Gauss sum problem is to evaluate $G(\q)$ with
the first examples going back to Gauss \cite{bG}*{Art.356}.

In the generality of a function on a finite set, it is difficult
to say very much useful except that the problem
can be divided into a magnitude and a phase.
We can consider the magnitude or norm, 
$\magg\q=\vert G(\q)\vert$,
and when $\magg\q\neq0$, define the phase 
$\arf(\q)\in\rmodz$
by
\[G(\q)=\magg\q\cdot \exp{\arf(\q)}\ .\]
Here are two general constructions. 
Given $\q_i$ on $T_i$, $i=1$, $2$, define the \emph{orthogonal\/} sum
$\q_1\perp\q_2\colon T_1\times T_2\to\rmodz$
by $(\q_1\perp\q_2)(x_1,x_2)=\q_1(x_1)+\q_2(x_2)$.
Check that

\[G(\q_1\perp\q_2)=G(\q_1)\cdot G(\q_2)\ .\leqno(\nextLabel{})\]

For the second construction, 
observe $\Z$ acts on the functions by $(a\cdot\q)(x)=a\cdot\q(x)$
for $a\in\Z$.
For $a=-1$,

\[G(-\q)=\conj{{G(\q)}}\ .\leqno(\nextLabel[xxB]{})\]

We review and extend some examples for which
the magnitude and phase can be calculated and apply these
results to some problems in algebra and topology.

Start with an easy example: $T$ is a finite group
and $\q$ is a homomorphism.
If $\q$ is trivial, the answer is immediate:
$\magg\q=\vert T\vert$ and $\arf(\q)=0$ 
since every term in the sum is $1$.
For the non{--}trivial case do the sum in two ways and
compare the answers: specifically, fix an element $c\in T$ with $\q(c)\neq0$ 
and observe
\[\sum_{t\in T} \exp{\q(t)}=\sum_{t\in T} \exp{\q(t+c)}\]
because we are summing over the same set and
\[\sum_{t\in T} \exp{\q(t+c)}=\exp{\q(c)}
\sum_{t\in T} \exp{\q(t)}\]
since $\q$ is a homomorphism. 
Since $\exp{\q(c)}\neq 1$, 
$\magg\q=0$.

Our main interest is the case where $T$ is a finite abelian
group and $\q$ has some sort of ``quadratic'' property, which we now
review.
A function $\q\colon T\to \qmodz$ is called a 
{\sl quadratic function\/}
provided
\begin{itemize}
\item[{(\nextLabel{})}] $\q(a x)=a^2\q(x)$ for all integers $a$ 
and all $x\in T$
\item[{(\nextLabel[bilinearEq]{})}] $\b(x,y)=\q(x+y)-\q(x)-\q(y)
\colon T\times T\to\qmodz$ defines a bilinear form.
\end{itemize}
\noindent
We call $\q$ a {\sl quadratic enhancement\/} of $b$ and $b$ the
{\sl associated bilinear form\/} to $\q$.
If $\q$ just satisfies \bilinearEq, we say $\q$ is an 
{\sl enhancement\/} of $b$.
Condition \bilinearEq\ is equivalent to the condition

\begin{align*}
(\nextLabel{})\quad\q(x_1+x_2+x_3)=&
\q(x_1+x_2)+\q(x_1+x_3)+\q(x_2+x_3)-\\&
\bigl(\q(x_1)+\q(x_2)+\q(x_3)\bigr)
\end{align*}
for all $x_i\in T$.
Any function $\q\colon T\to\rmodz$
satisfying \bilinearEq\ will be called {\sl $2${--}linear}.
The \emph{associated bilinear form\/}, $\abf{\q}$,
is defined by \bilinearEq: 
note that the bilinear form is symmetric.
The set of functions from $T$ to $\rmodz$ is a group
and the set of $2${--}linear functions is a subgroup.
The orthogonal sum of $2${--}linear functions is $2${--}linear.

Any symmetric bilinear form has enhancements.
Given $\b\colon T\times T\to\qmodz$ define
$\wq \b$ to be $T\oplus\qmodz$ with group structure
$(t_1,r_1)+(t_2,r_2)=\bigl(t_1+t_2, r_1+r_2-\b(t_1,t_2)\bigr)$.
It is straightforward to check that $\wq b$ is an abelian group and
that $\iota\colon\qmodz\to\wq b$ defined by $\iota(r)=(0,r)$
is an injective homomorphism.
Since $\qmodz$ is divisible, $\iota$ has a section,
$\Psi\colon\wq b\to\qmodz$.
Define $\q(t)=\Psi(t,0)$ and check that $\q$ is an enhancement of $\b$.
Conversely, given an enhancement of $\b$, define
$\Psi(t,m)=\q(t)+m$ and check that it is a section, so the set
of enhancements of $\b$ corresponds bijectively 
to the set of sections of $\iota$.
Let $\enh \b$ denote the space of sections of $\iota$
or equivalently the space of enhancements of $\b$.

The set $\enh\b$ is acted on by the group
of homomorphisms, $T^\ast=\dual{T}$: 
given one enhancement of $\b$, say $\q$,
 and $h\in T^\ast$ define $\q_h$ by $\q_h(t)=\q(t)+h(t)$
and check that $\q_h$ is an enhancement of $\b$.
If $\q_1$ and $\q_2$ are enhancements of $\b$, check that
$\q_1-\q_2\in T^\ast$.
Pick an element $\q\in\enh\b$ and define
$T^\ast\to\enh\b$ using the action on $\q$:
this function is a bijection.

A  homomorphism, $h\colon T_1\to T_2$ is an 
{\sl isometry\/} between
two $2${--}linear functions provided  $\q_2\bigl(h(x)\bigr)=\q_1(x)$:
if $h$ is an isomorphism we say $\q_1$ and $\q_2$ 
are {\sl isometric\/}, written $\q_1\cong\q_2$.
The map $h$ is an isometry in the usual sense between 
the associated bilinear forms.
As an example, if $a\in Z$ is relatively prime to $\vert T\vert$ then
multiplication by $a$ on $T$ gives an isometry between $a^2\cdot\q$
and $\q$ if $\q$ is quadratic.

Any $2${--}linear function satisfies $\q(0)=0$ and
$b(x,x)=\q(x)+\q(-x)$ (compute $b(x,-x)$
using bilinearity).
Induction shows that
\[\q(a x)=\frac{a^2+a}{2}\cdot\q(x)+
\frac{a^2-a}{2}\cdot\q(-x)\ .\leqno(\nextLabel[xxxB]{})\]
Furthermore, $\qh\q\colon T\to\qmodz$
defined by $\qh\q(x)=\q(x)-\q(-x)$ is a homomorphism
and $\q$ is quadratic if and only if $\qh\q$ is identically $0$.
Check that $\qh{\q}$ vanishes on elements in $T$ of order $2$
and that $\qh{\q}$ behaves well under our two constructions:
$\qh{(\q_1\perp\q_2)}=\qh\q_1+\qh\q_2$ and
$\qh{(a\cdot\q)}=a\cdot\qh{\q}$.

It follows that,  if $ax=0$ then
$a\q(x)\pm a\q(-x)=0$.
Hence $2a\q(x)=0$ and any $2${--}linear function 
automatically lands in $\qmodz$.
If $a$ is odd it further follows that $a^2\q(x)=0$ and hence
$a\q(x)=0$.
It follows that $\q=\perp_p\q_p$,
where $\q_p$ denotes $\q$ restricted to the
$p${--}torsion subgroup and $\q_p$ takes values
in 
$\Z[\smfrac1p2]/\Z\subset\qmodz$.
Finally, it follows that the group of $2${--}linear functions is finite.

\medskip
{The action of $T^\ast$ on enhancements has the following effect.
Check that $\qh{(\q_h)}=\qh{\q}+2h$ so if $\q$ is quadratic,
$\q_h$ is also quadratic if and only if $h$ 
takes values in $\cy2\subset\qmodz$.
Since $\qh{\q}$ vanishes on elements of order $2$, there exits a homomorphism
such that $\qh{\q}=2h$ and then $\q_{-h}$ is quadratic.
If $\q$ is a quadratic enhancement of $\b$ then $\q_h$ is also quadratic 
if and only if $h$ takes values in $\cy2\subset\qmodz$.
}

\medskip
If $K\subset T$ is a subgroup, define 
$K^\perp=\{t\in T\ |\ \abf{\q}(t,k)=0\ \forall\ k\in K\}$.
We call $\q$ and $\abf{\q}$  {\sl non{--}singular\/} 
provided $T^\perp=\{0\}$.
Note $(T_1\perp T_2)^\perp=T_1^\perp\oplus T_2^\perp$.
If $K$ satisfies $\q\vert_K=0$, then $K\subset K^\perp$ and
$\q$ induces a well{--}defined function
\[\ir\q K\colon K^\perp/K\to\qmodz\ .\]
Check that
$\bigl( K^\perp/K\bigr)^\perp=(T^\perp)/(T^\perp\cap K)$
so $\ir\q K$ is tame (resp. non{--}singular) if $\q$ is.

Any bilinear form $\b$ defines an adjoint homomorphism
$\ad(\b)\colon T\to\dual{T}$ and $T^\perp$ 
is the kernel of $\ad(\b)$.
This means that non{--}singular forms have the property
that for any $h\in T^\ast$ there exists a unique $c\in T$
such that $h(t)=\b(t,c)$ for all $t\in T$.
In the singular case, given any $h\in T^\ast$
with $h\vert_{T^\perp}$ trivial,
there exist $c\in T$, not unique, such that
$h(t)=\b(t,c)$ for all $t\in T$.
\begin{DefS}{\nextLabel[subx]{Definition}}
For $x\in T/T^\perp$ we use the notation $\q_x$ to denote the enhancement
$\q_x(t)=\q(t)+\b(x_1,t)$ for all $t\in T$, where $x_1\in T$ maps
to $x$.
\end{DefS}

Check that $\q_x$ is independent of the choice of $x_1$ as is $\q(x_1)$,
which will be denoted $\q(x)$.

\begin{ThmS}{\nextLabel[subxpsi]{Lemma}}
$\q_x(t) = \q(t+x_1) - \q(x_1)$ so 
$\arf(\q_x)=\arf(\q)-\q(x)$. 
\end{ThmS}

For any $2${--}linear function, $\q$ restricted to $T^\perp$ 
is a homomorphism: denote it by $\qs{\q}$.
We say $\q$ is {\sl tame\/} provided $\qs{\q}$ 
is trivial.
If $\q$ is tame then it induces a non{--}singular form on $T/T^\perp$,
which we denote by $\red\q$
The orthogonal sum of two tame $2${--}linear functions is tame.
If $\q$ is tame, so is $a\cdot\q$ for any $a\in\Z$ which is relatively
prime to $\vert T\vert$.
If $\q$ is quadratic $\qs{\q}$ takes values in $\cy2\subset\qmodz$
and whenever $2x\in T^\perp$, $\qs{\q}(2x)=0$.

Suppose $\q$ is tame. Then $\q_h$ is tame if and only if 
$h$ vanishes on $T^\perp$.

\begin{DefS}{\nextLabel[tf]{Remarks}}
Since $\qmodz$ is injective 
there is always an extension $h$ of $\q^s$ to all of $T$ and for 
any such $h$, $\q_{-h}$ is tame.
If $\q$ is tame the function $x\in T/T^\perp\mapsto \q_x$
defines a bijection between $T/T^\perp$ and the tame enhancements
of $\b$. 
\end{DefS}

\begin{DefS}{\nextLabel[qf]{Remarks}}
For any abelian group $T$ and $n\in\Z$,
let $\preT n T=\bigl\{ x\in T\ \vert\ n x=0\bigr\}$ and
let $\dualT n T=\bigl\{ x\in T\ \vert\ x=n y \bigr\}$.
If $\q$ is tame and quadratic the function 
$x\in\preT2{(T/T^\perp)}\mapsto \q_x$
defines a bijection between $\preT2{(T/T^\perp)}$ 
and the tame quadratic enhancements
of $\b$.
\end{DefS}

\begin{ThmS}{\nextLabel{Theorem}}
Every symmetric bilinear form $b$ has a tame quadratic enhancement.
\end{ThmS}
\begin{proof}
We have seen $\b$ has enhancements.
Using (\tfLink) construct a tame one $\q_t$.
Pass to $\red{\q_t}$ and use (\qfLink) to construct 
a quadratic $\hat\q\colon T/T^\perp\to\qmodz$
and check that the composition
$\vrule width0pt height12pt depth 0pt
\q\colon T\to T/T^\perp\ \xrightarrow{\ \hat\q\ }\ \qmodz$
is a tame quadratic enhancement of $\b$.
\end{proof}

\bigskip

For subgroups $K\subset T$ with $\q\vert_K=0$,
the Gauss sums for $\q$ and $\ir\q K$ are related.
This has been noticed before,
\cite{bB}, \cite{bBM}, \cite{bJM}, \cite{bMH}, \cite{bT} and many others, 
but apparently only for non{--}singular quadratic functions.

\begin{ThmS}{\nextLabel[xxC]{Theorem}}
Let $\q\colon T\to\rmodz$ be $2${--}linear and suppose $\q\vert_K=0$
for some subgroup $K$.
Then
\[G(\q)=\vert K\vert\cdot G(\q|_{K^\perp/K})\ .\]
\end{ThmS}
We emphasize that neither tame nor quadratic is assumed.
\begin{proof}
The usual proof assuming non-singular still works.
Pick coset representatives $\alpha_i$ for $T/K^\perp$ and
$\gamma_j$ for $K^\perp/K$ and let $\kappa_k$ run over the
elements of $K$.
Then every element of $T$ can be written uniquely as
$\alpha_i+\gamma_j+\kappa_k$ and
$\q(\alpha_i+\gamma_j+\kappa_k)=
b(\alpha_i,\gamma_j)+b(\alpha_i,\kappa_k)+\q(\alpha_i)+\q(\gamma_j)$.
If we fix $\alpha_i$ and $\gamma_j$ and sum over $\kappa_k$
then the only term which depends on $\kappa_k$ is $b(\alpha_i,\kappa_k)$
which is a homomorphism.
Hence the sum is $0$ if $b(\alpha_i,\kappa_k)$ is non{--}zero for some
$\kappa_k$.
But by definition of $K^\perp$, $b(\alpha_i,\kappa_k)$ is non{--}trivial
except for the $\alpha_i$ in the $0$ coset, say $\alpha_0$, and in this case
the sum is just $\vert K\vert\cdot\exp{\q(\alpha_0+\gamma_j)}$.
\end{proof}

A standard observation permits us to evaluate the magnitude of $G(\q)$.
\begin{ThmS}{Theorem \nextLabel[xxD]{Theorem}}
If $\q$ is $2${--}linear then $\magg{\q}=0$ if $\q$ is not tame and 
$\magg{\q}=
\sqrt{\vert T^\perp\vert\cdot\vert T\vert\hskip1pt}$
 if $\q$ is tame.
\end{ThmS}
\begin{proof}
If $\q$ is not tame then there is a $c\in T^\perp$ such that
$\q(c)\neq0$ and $\q(t+c)=\q(t)+\q(c)$ for all $t\in T$.
As in the homomorphism case it follows that $\magg{\q}=0$.
If $\q$ is tame, 
check that the evident inclusion
$\Delta\colon T\subset T\perp T/T^\perp$ 
with 2{--}linear function $\q\perp-\red\q$ satisfies 
$\bigl(\Delta(T)\bigr)^T=\Delta(T)$ and $\q\perp-\red\q$ 
restricted to $\Delta(T)$ vanishes so
$G(\q\perp-\red\q)=\magg{T}$.
But $G(\q\perp-\red\q)=G(\q)\cdot G(-\red\q)$;
$G(\q)=\magg{T^\perp}\cdot G(\red\q)$;
and $G(-\red\q)=\Bar{G(\red\q)}$.
\end{proof}

We can evaluate $\arf(\q_h)$ in terms of
$\arf(\q)$ if both are tame.
\begin{ThmS}{\nextLabel[xxE]{Theorem}}
Suppose $\q$ is $2${--}linear and tame and that $h\colon T\to\qmodz$
is a homomorphism.
Then $\q_h$ is tame if and only if $h(T^\perp)=0$.
If $q_h$ is tame
let $c\in T/T^\perp$ be the unique element corresponding to $h$,
so $\q_h=\q_c$.
The value of $\q(c_1)\in\qmodz$ is the same for all elements reducing to
$c$, so let us denote that common value by $\q(c)$.
Then
\[\arf(\q_h)=\arf(\q)-\q(c) .\]
\end{ThmS}
\begin{proof}
Since $\q_h=\q_c$, the displayed formula holds by \subxpsiName\ \subxpsiLink.
\end{proof}

\gdef\mymu{\Delta_\arf}
\begin{DefS}{\nextLabel[dmd]{Remark}}
We check that a potential source of new functions in fact yields
nothing new.
Let $\q$ be a tame enhancement and define a new function
$$\mymu\colon T/T^\perp\to\rmodz$$
by $\mymu(x)=\arf(\q)-\arf(\q_x)$ for all $x\in T/T^\perp$.
Then (\xxELink) says $\mymu=\red\q$.
If $\q$ is additionally quadratic then 
$\preT2(T/T^\perp)$ acts simply{--}transitively
on the tame quadratic enhancements.
In example (\xxAfLink) below  we encounter 
the construction $\mymu$ and it is nice to be able to identify it.
\end{DefS}

\bigskip

It follows from a beautiful argument due to 
Frank Connolly, \cite{bC}*{p.393},
that, in the quadratic case, 
$\q\perp\q\perp\q\perp\q$ is isometric to
$-\q\perp-\q\perp-\q\perp-\q$.
It follows that for any tame, quadratic $2${--}linear function, 
$\arf(\q)\in\cy8\subset\qmodz$ and 
it then follows from (\xxELink) that $\arf(\q)\in\qmodz$ for any
tame $2${--}linear function.

Connolly's result can be made more precise: let
$T_p$ denote the $p${--}torsion subgroup of $T$
and let $\q_p=\q\vert_{T_p}$.
As usual, it suffices to understand $\q_p$.

\begin{ThmS}{\nextLabel[fxc]{Theorem}}
Assume $\q_p$ is quadratic $2${--}linear.
If $p\equiv 1\bmod 4$, $\q_p\cong-\q_p$;
if $p\equiv 3\bmod 4$, $\q_p\perp\q_p\cong-\q_p\perp-\q_p$;
if $p-2$, $\q_2\perp\q_2\perp\q_2\perp\q_2
\cong -\q_2\perp-\q_2\perp-\q_2\perp-\q_2$.
\end{ThmS}

Here is a version of Connolly's argument.
Choose $k$ so that $p^k$ is an exponent for $T_p$
and hence $\q_p$ takes values in $\cy{p^k}\subset\qmodz$,
$p$ odd or in $\cy{2^{k+1}}\subset\qmodz$.
Find integers $x_i$ such that
$x_1^2+x_2^2+x_3^2+x_4^2\equiv -1\bmod p^\ell$ with
$\ell\geq k$ $p$ odd or $\ell> k+1$ if $p=2$.
Additionally assume $x_3=x_4=0$ if $p$ is odd and further assume $x_2=0$
if $p\equiv1 \bmod 4$.
Assuming this done,
Connolly's matrix $M=\left[%
\hskip6pt\vrule height 24pt depth 24pt width 0pt
\begin{matrix}
x_1&x_2&x_3&x_4\cr
\m x_2&x_1&\m x_4& x_3\cr
x_3&\m x_4&\m x_1&x_2\cr
x_4&x_3&\m x_2&\m x_1\cr
\end{matrix}\right]$
corresponds to an evident linear map from 
$T_p\oplus T_p\oplus T_p\oplus T_p$
to itself.
Check $M\cdot M^T=(x_1^2+x_2^2+x_3^2+x_4^2)I$ where $I$ is the 4 by 4 identity matrix. 
Hence $M$
defines a bijection which
can be checked to give an isometry between
$\q_p\perp\q_p\perp\q_p\perp\q_p$ and 
$-\q_p\perp-\q_p\perp-\q_p\perp-\q_p$.
If $x_3=x_4=0$  or $x_2=x_3=x_4=0$ the evident
square submatrix
gives the required isometry.

Connolly produced the $x_i$ by an appeal to Lagrange's
four{--}squares theorem.
Here is a different approach which uses only quadratic reciprocity.
If we could solve
\[x_1^2+x_2^2+x_3^2+x_4^2\equiv -1\bmod p^\ell\leqno(\nextLabel[xxF]{})\]
for all large $\ell$ we would be done.

If $(x_1, x_2, x_3, x_4)$ satisfies (\xxF) for $p$
odd and $\ell\geq 1$,
then we can find $a$ such that
$(x_1+a p^\ell)^2+x_2^2+x_3^2+x_4^2\equiv -1
\bmod p^{\ell+1}$.
If $p=2$ and if
$(x_1, x_2, x_3, x_4)$ satisfies (\xxF) for $p=2$
and $\ell\geq 3$,
then we can find $a$ so that
$(x_1+a2^{\ell-1})^2+x_2^2+x_3^2+x_4^2\equiv -1
\bmod 2^{\ell+1}$.
This construction is just Hensel's lemma for this simple case.

Note $(1,1,1,2)$ satisfies (\xxF) for $p=2$, $\ell=3$.
If $p\equiv1\bmod4$, then $-1$ is a quadratic residue
so we can solve $x_1^2\equiv-1\bmod p$.
If $p\equiv3\bmod 4$,
$x_1^2+x_2^2\equiv-1\bmod p$ has solutions
since there are $\frac{p+1}2$ distinct values for $-1-x_2^2$ and 
only $\frac{p-1}2$ non{--}residues.\qed

Brown \cite{bB} studied the case in which $T$ is a $\cy2$ vector space.
Any $2${--}linear function on a $\cy2$ vector space is quadratic and
Brown's functions were assumed non{--}singular, although tame would 
have sufficed for many of his results.
For example, Brown gave a different argument for 
$\arf(\q)\in\cy8$.
An element $\c\in T$ is {\sl characteristic\/} provided
$b(x,x)=b(\c,x)$ for all $x\in T$.
Brown's argument is to observe that characteristic elements
exist and then $\q_{\c}=-\q$.
By (\xxELink), if $\q$ is tame, $\arf(\q)=\q(\c)+\arf(\q_\c)\in\qmodz$,
so $2\arf(\q)=\q(\c)\in\cy4\subset\qmodz$.
\
For later use, recall the Gauss sum formula of Dirichlet
\cite{bD} that we require.

\begin{ThmS}{\nextLabel[dgsf]{Theorem}}
If $m>0$ then 
\[\sum_{s=0}^{m-1}
e^{2\pi i s^2/m}=\begin{cases}%
(1+i)\sqrt{m}& m\equiv 0\ \MOD 4\cr
\sqrt{m}& m\equiv 1\ \MOD 4\cr
0& m\equiv 2\ \MOD 4\cr
i\sqrt{m}& m\equiv 3\ \MOD 4\cr
\end{cases}\]
\end{ThmS}
Some elementary Galois theory allows us to extend (\xxxBLink).
Let $\q\colon T\to\qmodz$ be a quadratic function on a finite $p${--}group
of order $p^r$ and let $a\in\Z$ be prime to $p$.
As we saw above, if $a_1\equiv a_2\cdot s^2$ \MOD{p^{r+1}}, 
then $\arf(a_1\cdot\q)=\arf(a_2\cdot\q)$.
For odd $p$ the integers \MOD{p^{r+1}} divide into two sets, the quadratic
residues and the quadratic non{--}residues and the class to which an
integer belongs can be determined by examining it \MOD{p}.
For $p=2$ life is slightly more complicated: given odd integers 
$a_1$ and $a_2$, the equation $a_1\equiv a_2\cdot s^2$ 
\MOD{2^r} can be solved for all $r$ if and only if $a_1\equiv a_2$
\MOD{8}.

We can now work out the formula for $\arf(a\cdot\q)$ when $a$ is prime
to $p$. 
(For $p^i\cdot\q$ see (\oDqLink)  and (\oDrLink)).

\begin{ThmS}{\nextLabel[pof]{Theorem}}
Let $T$ be a $p$ group with a tame quadratic enhancement $\q$
and let $a$ be an integer prime to $p$.
Let $\vert T/T^\perp\vert=p^e$.
For $p$ odd, let $\qrN p(a)=0$ if $a$ is a quadratic residue \MOD{p} 
and $\qrN p(a)=1$ if $a$ is not.
Let $\qrN 2(a)=0$ if $a\equiv\pm1$ \MOD{8} and $\qrN 2(a)=1$ if
$a\equiv\pm3$ \MOD{8}.
Then for $p$ odd
\[\arf(a\cdot\q)=\arf(\q)+\frac {\qrN p(a)\cdot e}2\leqno(\pof_p)\]
and for $p=2$
\[\arf(a\cdot\q)=a\cdot\arf(\q)+\frac {\qrN 2(a)\cdot e}{2}\leqno(\pof_2)\]
\end{ThmS}
\begin{proof}
It suffices to do the non{--}singular case since 
we can work with $\red\q$.
Let $\zeta$ be a primitive $p^r$ root of unity where $p^r$ annihilates 
$\q(t)$ for all $t\in T$.
Let $\omega=\exp{\frac18}$.

Using \fxcLink\ we see 
\[G(\q)=p^{\frac e2}\omega^\sigma\hskip20pt \text{where } \arf(\q)=\frac{\sigma}{8}. 
\leqno(\nextLabel[gsf]{})\]

Now the left hand side of (\gsf) lies in $\Q[\zeta]$ 
and hence so must the right.
The integers relatively prime to $p$ map onto the Galois group of 
$\Q[\zeta]$ over $\Q$.
The map sends the integer $a$ to $\zeta^a$ so we need to study the effect of
the Galois automorphism $\gal{a}\colon \zeta\mapsto\zeta^a$ on $G(\q)$.

For $p$ odd, $\Q[\zeta]\cap\Q[\omega]=\Q$, so
$\pm1$ are the only powers of $\omega$ in $\Q[\zeta]$.
For $p\equiv 1$ \MOD4 Dirichlet (\dgsfLink) shows $\sqrt{p}\in\Q[\zeta]$
and one can check that $\gnqr(\sqrt{p})=-\sqrt{p}$ for some $\nqr$ in the Galois group.

If $e$ is even, then the right hand side of (\gsfLink)
must be an integer, so $\arf(\nqr\cdot\q)=\arf(\q)$.
If $e$ is odd, then the right hand side of (\gsfLink)
must be $\sqrt{p}$ times an integer, 
so $\arf(\nqr\cdot\q)=\arf(\q)+\frac12$.
Additionally $\arf(\q)$ is either $0$ or $\frac12$.
Note both cases are covered by the formula $(\pofLink)_p$.

For $p\equiv 3$ \MOD{4} Dirichlet (\dgsfLink) shows $i\sqrt{p}\in\Q[\zeta]$
and $\gnqr(\sqrt{p}i)=-\sqrt{p}i$.
If $e$ is even, then the right hand side of (\gsfLink)
must again be an integer, so $\arf(\nqr\cdot\q)=\arf(\q)$.
If $e$ is odd, then the right hand side of (\gsfLink)
must be $\sqrt{p}i$ times an integer, 
so $\arf(\nqr\cdot\q)=\arf(\q)+\frac12$ again and
$(\pofLink)_p$ holds in this case too.
Additionally $p\equiv3$ \MOD{4}, $\arf(\q)=\pm\frac{1}{4}$ if $e$ is odd and
$0$ or $\frac12$ if $e$ is even.

For $p=2$, there are three multiplications to be worked out.
Which class an integer $a$ belongs to can be determined by
reducing $a$ \MOD{8}.
The numbers are $-1$ and $\pm3$ \MOD{8}.
The affect of $-1$ we know: $\arf(-\q)=-\arf(\q)$.
The wrinkle when we multiply by $\pm3$ is that these
Galois actions send $\sqrt{2}$ to $-\sqrt{2}$.
If $e$ is even then this does not matter and
$\arf(\nqr\cdot\q)=\nqr\cdot\arf(\q)$
When $e$ is odd we get
$\arf(\nqr\cdot\q)=\nqr\cdot\arf(\q)+\frac12$.
So for $p=2$ $(\gsfLink)$ holds.
\end{proof}

\section{Some algebra applications.}\resetSectionLabel

\begin{ThmS}{\nextLabel[xxAa]{Theorem}}
Given a symmetric bilinear form $B$ 
on a rational vector space, $V$,
define $Q\colon V\to \Q$ by 
$Q(v)=\frac{\scriptstyle B(v,v)}{\scriptstyle2}$.
Call a lattice {\sl integral\/} if $B(v_1,v_2)\in\Z$ for all
$v_1$, $v_2\in L$.
Pick a lattice $L\subset V$ such that
$Q(x)$ is integral for all $x\in L$.
Define $L^{\#}=\{v\in V\ \vert\ B(v,\ell)\in\Z\ \forall\ \ell\in L\}$
and check that $L$ integral implies $L\subset L^{\#}$.
If $B$ is non{--}singular, check that $L^{\#}/L$ is finite.
Check that $Q$ induces a non{--}singular quadratic function 
$\q_L\colon L^{\#}/L\to\qmodz$
which enhances the symmetric bilinear form
$\b_L\colon L^\#/L \times L^\#/L\to\qmodz$ induced by $B$.
The Milgram Gauss Sum Formula \cite{bJM} says
$$\arf(\q_L)=\frac{\sig(B)}8\ ,\leqno(\xxAa)$$
where $\sig(B)$ denotes the signature of $B$.
\end{ThmS}
\medskip

We prove the formula using some 
straightforward manipulations and Dirichlet's Gauss sum formula.
The basic outline, except for the appeal to Dirichlet,
is in Milnor and Husemoller \cite{bMH} 
who attribute it to Knebusch.
A proof for $\det B$ odd was given earlier by Blij \cite{bBl}.

\begin{proof}
Call a lattice $L$ {\sl acceptable\/} if $Q(x)\in\Z$ for all $x\in L$.
If $L_1$ and $L_2$ are acceptable lattices, so is $L_1\cap L_2$.
To show all acceptable lattices give the same answer, it suffices 
to show $\arf(\q_{L_2})=\arf(\q_{L_1})$ under the
additional assumption that $L_1\subset L_2$, and hence
$L_1\subset L_2\subset L_2^{\#}\subset L_1^{\#}$.
Let $T=L_1^{\#}/L_1$ and apply (\xxCLink) to
$K=L_2/L_1\subset T$.
Check $K^\perp=L_2^{\#}/L_1$ so $K^\perp/K=L_2^{\#}/L_2$.

Over $\Q$, $B$ can be diagonalized and there are acceptable diagonal
lattices, so it suffices to show
$\arf(\rof{2m})=\frac18$ if $m>0$.
Now $G(\rof{2m})=\sum_{s=0}^{2m-1} \exp{\smfrac{s^2}{4m}{2}}$.
Dirichlet (\dgsfLink) in case $4m$ says
$\sum_{s=0}^{4m-1} \exp{\smfrac{s^2}{4m}{2}}=(1+i)\sqrt{4m}$.
Since $(s+2m)^2=s^2\ {\rm  mod}\ 4m$, 
$(1+i)\sqrt{4m}=2\cdot G(\rof{2m})$ and
the result follows.
\end{proof}

\medskip

\begin{ThmS}{\nextLabel[bmf]{Theorem}}
Suppose as above that $B$ is a non{--}singular bilinear form on a rational
vector space $V$ and $L$ is an integral lattice with $L^\#$ 
defined as in (\xxAaLink).
Check $B$ still induces a symmetric, bilinear
form $\b_L$ on $L^\#/L$ which is still non{--}singular. 
To apply the Milgram Gauss sum formula to $L^\#/L$ we need additionally
that $B(v,v)\in2\Z$ for all $v\in L$.
If $L$ does not satisfy this condition we can proceed as follows.
There exists a {\sl characteristic element \/}$\nu\in L$: i.e. 
$B(v,v)\equiv B(\nu, v)$ \MOD2 for all $v\in L$.
The element $\nu$ is not unique but pick one.
Then set $Q(x)=\frac{B(x,x)-B(\nu, x)}{2}$ for all $x\in L^\#$.
Check that $Q$ induces a quadratic function
$\q_\nu\colon L^\#/L\to\qmodz$ and that $\q_L$ 
enhances $\b_L$.
Then
\[\arf(\q_\nu)=\frac{\sig(B)}{8}-\frac{B(\nu,\nu)}{8}\in\qmodz\ .
\leqno(\bmf)\]
\end{ThmS}
\begin{proof}
Let $L_1=2L\subset L$.
Check that the function $Q$ defined above induces a function
$\q_{L_1}\colon L_1^\#/L_1\to\qmodz$ which is 
an enhancement of $\b_{L_1}$.
The function $\q_{L_1}$ is almost certainly not quadratic, but
just as in the proof of (\xxAaLink) prove that 
$\arf(\q_\nu)=\arf(\q_{L_1})$.
Note $\q_{L_1}+B( \_,\nu/2)=\q$ on $L_1^\#/L_1$ where $\q$ is
the quadratic function induced by $\frac{B(x,x)}{2}$ on $L_1^\#/L_1$.
Use (\xxAaLink)  to calculate $\arf(\q)=\frac{\sig(B)}{8}$ and
use (\xxELink) to deduce that
$\arf(\q_\nu)=\arf(\q)-\q(\nu/2)=
\arf(\q)-\frac{B(\nu,\nu)}{8}$.
\end{proof}

\begin{DefS}{\nextLabel[xxAb]{Remarks}}
If $\det B$ is odd, then $\nu$ is unique up to sums with elements
of the form $2x$.
A classical argument shows 
$B(\nu,\nu)\equiv B(\nu+2x,\nu+2x)$ \MOD8 and one checks that
$\q_\nu$ does not depend on the choice of $\nu$ either.
For $\det B$ odd (\bmfLink) is a result of Blij \cite{bBl}.
The general case appears in Brumfiel\&Morgan \cite{bBM}.

If $\det B$ is even however, there are different choices for $\nu$
which give different enhancements and different values of $B(\nu,\nu)$.
Indeed, it is a theorem of Brumfiel\&Morgan \cite{bBM}
and Wall \cite{bWa} that
any quadratic enhancement of a non{--}singular symmetric bilinear form
can be obtained as $\q_L$ for an appropriate $B$ and $\nu$.
\end{DefS}

\bigskip
We next turn to Brown \cite{bB} for other ways to obtain quadratic
functions.
Given a symmetric bilinear form, 
$B\colon V\times V\to\Z$, define 
$$\q_B\colon V\otimes\cy2\to\cy4\subset\qmodz$$
by $\displaystyle\q_B(x)=\frac{B(x,x)}4$.
This function is $2${--}linear, quadratic and the associated bilinear form is
\hbox{$B_2\colon (V\otimes\cy2)\times(V\otimes\cy2)\to\cy2$}
obtained by reducing $B$ \MOD{2}.
There is an obvious generalization: let
$V_m=V\otimes\cy{m}$ and define
$$\q_{B,m}\colon V_m\to\cy{2m}\subset\qmodz$$
with $\displaystyle\q_{B,m}(x)=\frac{B(x,x)}{2m}$.
In order for $\q_{B,m}$ to be defined on $V_m$ it is
necessary and sufficient that $m$ be even.
This quadratic enhancement, even for $m=2$, need not be tame.
It is non{--}singular if and only if $\det B$ is relatively prime to $m$.
A further generalization is to recall that for any 
quadratic $\q\colon T\to \qmodz$ the function
$B\otimes\q\colon V\otimes T\to\qmodz$ defined by
$(B\otimes\q)\bigl(v\otimes x\bigr)=B(v,v)\q(x)$ is also quadratic
\cite{bMH}*{p.111}.
(One should also check that it really is defined.
Note that the formula for $B=\rof1$ implies that $\q$ is quadratic
so the formula does not work for non{--}quadratic enhancements.)
If $m$ is even, let $\uf m\colon\cy{m}\to\qmodz$ be defined by
$\uf m(1)=\frac{1}{2m}$. 
Then $B\otimes\uf m=\q_{B,m}$ so this generalizes Brown's 
construction.

\begin{ThmS}{\nextLabel[xxAcx]{Theorem}}
If $B$ is a symmetric form over $\Z$ with determinant $\pm1$
and if $\q$ is tame quadratic,
then \[\arf(B\otimes\q)= \sig(B)\cdot\arf(\q)\in\qmodz\ .
\leqno(\xxAcx)\]
\end{ThmS}
\begin{proof}As Brown remarks, (\xxCLink) shows that
$\arf(B\otimes\q)$ only depends on the Witt class of $B$.
Now the Witt ring of $\Z$ is infinite cyclic generated by the 
form $\rof{1}$.
But $\rof{1}\otimes\q=\q$.
\end{proof}

\begin{ThmS}{\nextLabel[ebC]{Theorem}}
Another application of Dirichlet's Gauss sum formula  as in (\xxAaLink)
shows that for positive even $m$, $\arf(\uf m)=\frac{1}{8}$,
so we get Brown's theorem, $\arf(\q_{B,2})=\frac{\sig(B)}{8}$.
\end{ThmS}

\begin{ThmS}{\nextLabel[ebD]{Theorem}}
If $m$ is odd, define $\uf{m}\colon\cy m\to\qmodz$ by
$\uf{m}(1)=\frac1m$.
If $m>0$ and $m\equiv1$ \MOD{4}, $\arf(\uf m)=0$ and if
$m>0$, $m\equiv3$ \MOD{4}, $\arf(\uf m)=\frac14$.
\end{ThmS}
\bigskip

In a bit, we will give another proof of (\xxAcxLink) that bypasses the
Witt ring calculation, but 
before that, we generalize a formula of Melvin and Kirby
\cite{bKM}*{p.522} for computing $\arf(\q_{B,2})$ and use it to compute
$\arf(B\otimes\q)$ in favorable situations.

Given $\q$, fix a positive integer $m$ so that all the values of $\q$
applied to elements of $T$ lie in $\cy m$.
When computing $\arf(B\otimes\q)$ any two integer matrices which are the
same \MOD{m} clearly yield the same result.
We can also split $\q$ into its $p${--}primary pieces and work with
one prime at a time since $(B\otimes\q)_p=B\otimes \q_p$,
so let $m=p^k$.
We say $B$ is $p^r${--}similar to $C$, written
$B\sim_{p^r} C$, if there is an integral entry matrix $M$ 
with $\det M$ prime to $p$ so that $C\equiv M^{tr}BM$ \MOD{p^r}.
We note that if $B\sim_{m} C$ then
$M$ gives an isometry between $B\otimes\q$ and $C\otimes\q$, so
$\arf(B\otimes\q)=\arf(C\otimes\q)$ or neither is tame.

Call any integral matrix of the form 
$H_{m_1,m_2}=\begin{pmatrix}2m_1&1\cr1&2m_2\cr\end{pmatrix}$
{\sl pseudo{--}hyperbolic}.
Call any form which is an orthogonal sum of rank one forms and
pseudo{--}hyperbolics {\sl reduced}.
A form which is an orthogonal sum of rank one forms will be called
{\sl diagonal}.

The usual Gram{--}Schmidt process can be applied \MOD{p^r}
to yield the following algorithm for finding a reduced matrix ${p^r}${--}similar
to $B$.
If there is a diagonal entry $\alpha$ in $B$ generated by $x$ with $\alpha$
relatively prime to $p$, the usual Gram{--}Schmidt process defines
a projection $p\colon V\otimes\Q\to V\otimes\Q$ 
by $p(v)=v-\frac{\b(x,v)}{\b(x,x)}x$ and then observes that $B$ is
isometric to $\rof{\alpha}\perp \IM(p)$.
The only denominators in the matrix $B^\prime$ for $\IM(p)$ 
are divisors of $\b(x,x)$ so we can find a positive integer $r\equiv 1$ 
\MOD{p^r} so that $r B^\prime$ is integral and
$\rof{\alpha}\perp (rB^\prime)$ is congruent \MOD{m} to $B$.
Note $B$ and $\rof{\alpha}\perp (r B^\prime)$ are $p^r${--}similar.
Continue until we get to a new form 
$\bar B^\prime=D_0\perp B_0$ where $D_0$ is diagonal,
$B_0$ has all diagonal entries divisible by $p$ and
$\bar B^\prime$ is $p^r${--}similar to $B$.

Next suppose that some entry in $B_0$ is prime to $p$.
This means there are basis elements $x$, $y$ such that
$B(x,y)$ is prime to $p$.
If $p$ is odd change the basis to $x+y$, $x-y$ and note that
there are now diagonal entries prime to $p$.

If $p=2$ this does not work, but we can change the basis
element $y$ to $a y$, $a$ odd so that in this new basis
$B(x,y)=1$.
One can now orthogonally split off a  pseudo{--}hyperbolic.
Since before we split off the pseudo{--}hyperbolic all squares
were even afterwards all diagonal entries will still be even.

We can continue these reductions until we have found a matrix
$\bar B_0=R_0\perp B_1$ with $\bar B_0$ $p^r${--}similar to $B$,
with $R_0$ reduced and with every entry of $B_1$ divisible by $p$.
Divide $B_1$ by $p$ and continue.
Eventually we will obtain an orthogonal sum

\[C=C_0\perp pC_1\perp\cdots \perp p^{w}C_{w}
\leqno(\nextLabel[CrB]{})\]
where $C$ is $p^r${--}similar to $B$ and each $C_i$ is reduced
with $\det C_i$ prime to $p$.
If we take $r$ so that $p^r$ is at least the largest power of $p$ dividing
$\det B$ then $w$ and the sizes of each of the $C_i$ are determined
by $B$ since the $p$ torsion in cokernel of $C$ determines the $w$
and the size of the $C_i$ and for this large an $r$ the $p$ torsion in the
cokernel of $C$ is isomorphic to the $p$ torsion in the cokernel of $B$.

\bigskip
For $\det A$ prime to $p$ define a \MOD 8 integer $\sig_p(A)$ by

\[\sig_p(A)\equiv\begin{cases}
\rank A &\text{if $p$ is odd}\cr
N_1-N_{-1}+3N_{3}-3N_{-3}&\text{if $p=2$}\cr
\end{cases}\hskip8pt\biggl\}\hskip10pt\MOD8\]
where $N_i$ is the number of diagonal entries congruent to $i$ \MOD8
in any reduced matrix $8${--}similar to $A$.
It is not obvious that $\sig_p(A)$ is well{--}defined for $p=2$ 
but this is checked in the proof of the next result.

\begin{ThmS}{\nextLabel[lsmx]{Theorem}}
For each $p$ and $\q\colon T\to\qmodz$ with 
$\q$ tame and quadratic, $T$ a $p${--}group with 
$\vert T/T^\perp\vert=p^e$,
and $B$ a symmetric integral form
with $\det B$ prime to $p$,
\[\arf(B\otimes\q)=\sig_p(B)\cdot\arf(\q)+
\frac{\qrN p(\det B)\cdot e}{2}\ .\]
\end{ThmS}
\begin{proof}
We may pass to $\red\q$ so without loss of generality
assume $\q$ in also non{--}singular.
Note $B\otimes\q$ is also non{--}singular.
We start by proving the formula under the additional assumptions
that $B$ is reduced and we have define $\sig_p(B)$ 
using $B$ as our reduced form if $p=2$.
The desired formula is additive for orthogonal sum so it suffices to
prove the result for $\rof{a}$ plus the pseudo{--}hyperbolics if $p=2$.

Since $a\cdot\q=\rof{a}\otimes\q$ the formula for rank one
forms is just a restatement of the formula in \pofLink.

Next let $H$ be a pseudo{--}hyperbolic.
Working \MOD8 one can show that
\[\rof{-1}\perp\begin{pmatrix}2m_1&1\cr1&2m_2\cr\end{pmatrix}\sim_8
\rof{2m_1-1}\perp\rof{2m_2-1}\perp\rof{1-2m_1-2m_2}\ .
\]
Compute that
\[\arf(H\otimes\uf{2})=0\cdot\arf(\uf{2})+
\frac {\qrN 2(\det H)\cdot e}{2}\ .\]
Now work \MOD{2^r}, $r>3$, and suppose 
$\rof{-1}\perp\begin{pmatrix}2m_1&1\cr1&2m_2\end{pmatrix}\sim_{2^r}
\rof{a_1}\perp\rof{a_2}\perp\rof{a_3}$.
Calculating $\arf(H\otimes\uf{2})$ using this decomposition shows
$1+a_1+a_2+a_3\equiv0$ \MOD8.
It then follows that, for any tame $\q$,
$$\arf(H\otimes\q)=\arf(\q)\bigl(1+a_1+a_2+a_3)+
\frac {\qrN 2(a_1\cdot a_2\cdot a_3)\cdot e}{2}=
\frac {\qrN 2(\det H)\cdot e}{2}\ .$$
But $\sig_2(H)=0$ by definition. Hence the formula holds
for pseudo{--}hyperbolics.

Note next that if $B\sim_{p^r}C$, then 
$\qrN p(\det B)=\qrN p(\det C)$, provided $r\geq3$ if $p=2$.

If $p$ is odd but $B$ is not necessarily reduced, chose a reduced
$C$ with $C$ $p^r${--}similar to $B$ and with $p^r\cdot\q(t)=0$
for all $t\in T$.
Then $\arf(B\otimes\q)=\arf(C\otimes\q)$ and for odd primes clearly
$\sig_p(B)=\sig_p(C)$ so the formula holds for odd primes.

Now suppose $C_1$ and $C_2$ are two matrices $8${--}similar to $B$.
Then $\arf(B\otimes\uf{4})=\arf(C_i\otimes\uf{4})=\frac{\sig_2(C_i)}{8}$,
so $\sig_2(C_1)\equiv\sig_2(C_2)$ \MOD 8 
and we let $\sig_2(B)$ denote this
common value.
For any $\q$ and $B$ as  in theorem, choose a reduced 
$C$ $2^r${--}similar to $B$ with $2^r\cdot\q(t)=0$ for all $t\in T$
and just as in the odd case $\arf(B\otimes\q)=\arf(C\otimes\q)$.
We have just checked that $\sig_2(B)=\sig_2(C)$ and
$\qrN2(\det B)=\qrN2(\det C)$.
\end{proof}
\medskip

We can now prove (\xxAcxLink) directly from \lsmxLink.
It suffices to prove (\xxAcxLink) a prime at a time.
For $p\equiv1$ \MOD4, $\arf(\q)$ is a multiple of $\frac12$ and
$\qrN p(-1)=0$.
Since $\rank B\equiv\sig(B)$ \MOD2, (\xxAcxLink) holds.
For $p\equiv3$ \MOD4, $\arf(\q)$ is a multiple of $\frac14$:
recall $\rank B\equiv\sig(B)$ \MOD4 if $\det B>0$ and
$\rank B\equiv\sig(B)+2$ \MOD4 if $\det B<0$. 
Since $\qrN p(-1)=1$ (\xxAcxLink) holds again.

For $p=2$ we adopt a different approach.
Recall that we calculated
$\arf(B\otimes\uf{4})=\frac{\sig_2(B)}{8}$.
But we can calculate $\arf(B\otimes\uf{4})$ directly using the Milgram
Gauss sum formula.
First note that $2V$ is an acceptable lattice in $V\otimes\Q$
and since $\det B=\pm1$, $(2V)^\#=\frac12V$.
Calculate that $V^\#/V\cong V\otimes\cy4$ and under this isomorphism
the enhancement on $V^\#/V$ is the one on $B\otimes\uf{4}$.
Since $\det B=\pm1$, $\qrN 2(\det B)=0$ so
$\sig(B)\equiv \sig_2(B)$ \MOD{8}.\qed

\medskip

Kirby and Melvin \cite{bKM} introduced this reduction 
for $p^r=4$ including the
removal of pseudo{--}hyperbolics by orthogonally summing with
$\rof{1}\perp\rof{-1}$. 
If there are any diagonal entries of the form $\rof{2}$ they observe
that $\q_{B,2}$ is not tame; otherwise it is.
Then they compute that 
$$\arf(\q_{B,2})=\frac{n_1-n_{3}+4\epsilon_H}{8}\in\qmodz$$
where $n_1$ is the number of diagonal entries congruent to $1$ \MOD{4};
 $n_{3}$ is the number of diagonal entries congruent to $3$ \MOD{4};
and $\epsilon_H$ is the number of pseudo{--}hyperbolics congruent to
$\begin{pmatrix}2&\pm1\cr\pm1&2\cr\end{pmatrix}$ \MOD{4}.

For a proof, let $C\sim_4B$ where
$C=R_0\perp 2R_1\perp \cdots$
and $B\otimes\uf{2}=R_0\otimes\uf{2}\perp R_1\otimes (2\cdot\uf{2})$.
Now $n_2$ is the number of diagonal summands of $R_1$.
If $n_2>0$ the second summand is not tame 
and so neither is $B\otimes\uf{2}$.
Kirby and Melvin note that if $n_2>0$ then we can choose a $C$
so that $n_2$ is any number between $1$ and $\rank{R_1}+\rank{R_2}$.
If $n_2=0$ then $B\otimes\uf{2}=R_0\otimes\uf{2}$ orthogonal sum a $0$
form.
In this case $n_1=N_1+N_{-3}$ and
$n_3=N_3+N_{-1}$, so we would have
$\arf(\q_{B,2})=\arf(R_0\otimes\uf{2})=
\frac{n_1-n_3+4(N_3-N_{-3})}{8}+\frac{\qrN 2(\det R_0)}{2}$
and one checks that 
$\qrN 2(\det R_0)+N_3+N_{-3}$ is $\epsilon_H$ \MOD2.

Even though we can compute $\sig_2$ by working \MOD4 it is worth
considering the \MOD8 calculation since for $\det B$ odd we need
only do the reduction until all the diagonal elements are even, we need
not actually split off the pseudo{--}hyperbolics to compute $\sig_2$.
As an extreme example, if $B$ is a symmetric matrix over $\Z$ with all
diagonal entries even and with $\det B$ odd, then $\sig_2(B)=0$.

\begin{DefS}{\nextLabel[symmMat]{Remark}}
Given a symmetric matrix over $B$ with even diagonal 
entries we can find matrices $S$ such that $B=S+S^{tr}$.
\end{DefS}
The function
$\q\colon V\otimes\cy2\to\qmodz$ given by
$\q(x)=\frac{x^{tr}Sx}{2}$ is a quadratic enhancement of $B$:
in fact $\q=B\otimes\uf{2}$.
If $\det B$ is odd, it follows that
$\arf(\q)=\arf(\q_{B,2})=\frac{\qrN2(\det B)}{2}$.
If $S$ is a Seifert matrix for a knot $S^{4k-3}$ in $S^{4k-1}$, then
$\arf(\q)$ is the Arf invariant of the knot, $\det B=\Delta(-1)$ is the value
of the Alexander polynomial of the knot evaluated at $-1$ and
$\arf(\q)=\frac{\qrN2(\Delta(-1))}{2}$ 
is a theorem of Levine \cite{bL}*{p.544}.
\bigskip
In \cite{bT} we employed a slightly different variant of the Gram{--}Schmidt
process to classify arbitrary non{--}singular quadratic
functions on a $p$ group:
Thm. 3.5 (p.268) asserts that any non{--}singular $\q_p$ on
a finite abelian $p$ group is an orthogonal sum 

\begin{ThmS}{\nextLabel[oDq]{Theorem}}
\[\q_p=R_1\otimes\uf{p}
\perp\cdots\perp R_{k}\otimes\uf{p^k}\leqno(\oDq)\]
where each $R_i$ is reduced and has $\det R_i$ prime to $p$.
Furthermore,
$\rank{R_i}$ depends only on $T$
and 
\[\arf(\q_p)=\sum_{s=1}^{k}\arf(R_s\otimes\uf{p^{s}})=
\sum_{s=1}^{k}\ \Biggl(\sig_p(R_s)+\frac{\qrN p(\det R_s)\cdot s}{2}\Biggr)\ .\]
\end{ThmS}
\medskip

\begin{DefS}{\nextLabel[Bpt]{Remark}}
Brown \cite{bB} defined a product on enhancements on
$\cy2$ vector spaces.
We see no hope for defining a product in general, but given
two quadratic enhancements on $\cy p$ vector spaces, say
$\q_1=R\otimes\uf{p}$ and $q_2=S\otimes\uf{p}$, define
$\q_1\bullet\q_2$ to be $(R\otimes S)\otimes\uf p$.
Check that this product is well defined.
It clearly is commutative and it is not hard to verify that if
$\arf(\q_i)=\frac{\alpha_i}{8}$ then 
$\arf(\q_1\bullet\q_2)=\frac{\alpha_1\cdot\alpha_2}{8}$.
When $p=2$ this is Brown's product and Brown's theorem.
\end{DefS}

\bigskip
Given a form $B$ and a quadratic enhancement $\q$ we can work out
$(B\otimes\q)_p$ by finding a matrix ${p^r}${--}similar to 
one as in (\CrBLink) and a
quadratic enhancement $\q_p$ decomposed as in (\oDqLink), it
is clear that $(B\otimes\q)_p$ is an orthogonal sum
of terms of the form $(C_i\otimes R_j)\otimes(p^{i}\cdot\uf{p^{j}})$.
Hence, to describe $B\otimes\q$ or even just $p^i\cdot\q_p$
it suffices to describe $p^i\cdot \uf{p^{r}}$.
Let $T=\cy{p^r}$ denote the domain of $p^i\cdot\uf{p^r}$.
This form has $T^\perp=p^{r-i}T$.
If $p$ is odd $p^i\cdot\uf{p^r}$ is tame and 
$\red{(p^i\cdot\uf{p^r})}=\uf{p^{r-i}}$.
(Or trivial if $i\geq r$).
If $p=2$ this still works {\sl except\/} if $r=i$:
the form $2^r\cdot\uf{2^r}$ is not tame, although
$H_{m,n}\otimes(2^r\cdot\uf{2^r})=
2^r\cdot(H_{m,n}\otimes\uf{2^r})$ is trivial, hence tame.
Hence $2^i\cdot\q_2$ is tame if and
only if $R_{i}$ has no diagonal summands.
If $p^i\cdot\q_p$ is tame, we have

\begin{ThmS}{\nextLabel[oDr]{Theorem}}
\[\red{(p^i\cdot\q_p)}=R_{i+1}\otimes\uf{p}
\perp\cdots\perp 
R_{k}\otimes\uf{p^{k-i}}\leqno(\oDr)\]
\end{ThmS}

\begin{DefS}{\nextLabel[xxCC]{Remark}}
Note for $p$ odd, each $\arf(R_i\otimes\uf{p^{k-i}})$ is
an invariant of the isotropy class of $\q_p$. 
The proof is complicated a bit by the fact that although
$\displaystyle\arf(p^i\cdot\q_p)=\sum_{s=1}^{k-i}
\arf(R_{i+s}\otimes\uf{p^{s}})$,
$\arf(R\otimes\uf{p^s})$ does not in general equal
$\arf(R\otimes\uf{p^{s-1}})$ since the $e$ in (\lsmxLink)
keeps changing but still we have that
each $\sig_p(R_i)$ and each $\qrN p(\det R_i)$ depend only on
the isotopy class of $\q_p$.
If we have a symmetric matrix $B$ with $\det B$ prime to $p$
then if $B$ is $p^r${--}similar to a $C$ as in (\CrBLink) with
a large enough $r$, then $B$ determines each
$\sig_p(C_i)$ and each $\qrN p(\det C_i)$.
From this information, $\arf\bigl((B\otimes\q)_p\bigr)$ can be worked out.
\end{DefS}

If $\det A$ is prime to $p$, then $\sig_p(A)\equiv\rank A$ \MOD2.
For all primes $p$ one can check that for matrices with $\det(C\otimes R)$
prime to $p$, 
$\sig_p(C\otimes R)=\sig_p(C)\cdot\sig_p(R)$ and
$\qrN p(\det (C\otimes R))=
\sig_p(C)\cdot \qrN p(R)+\qrN p(C)\cdot\sig_p(R)$
computed \MOD2.

It follows that
\begin{ThmS}{\nextLabel[gtp]{Theorem}}
\begin{align*}\arf\bigl((B\otimes\q)_p\bigr)=&
\sum_{i<j}\frac{\sig_p(C_i)\cdot\sig_p(R_j)}8+\\&
\frac{\bigl(\sig_p(C_i)\cdot \qrN p(R_j)+
\qrN p(C_i)\cdot\sig_p(R_j)\bigr)\cdot(j-i)}{2}
\end{align*}\end{ThmS}

The result can be packaged a little better.

Let $\sig_p^\prime(A)=\sig_p(A)+4\cdot\qrN p(\det A)$. 
Then
\begin{align*}\hbox to 0pt{\hss$(\gtp)^\prime$\hskip40pt} \arf\bigl((B\otimes\q)_p\bigr)=&
\sum_{\substack{{i<j}\\{i\not\equiv j \ {\rm mod\ }2}}}
\frac{\sig_p(C_i)\cdot\sig_p(R_j)}8+\\&
\sum_{\substack{{i<j}\\{i\equiv j \ {\rm mod\ }2}}}
\frac{\sig_p^\prime(C_i)\cdot\sig_p^\prime(R_j)}8
\end{align*}

For $p=2$ the argument for invariance fails, as does the result: 
$\uf{4}\perp\uf{2}$ is isometric to 
$\rof{3}\otimes\uf{4}\perp\rof{3}\otimes\uf{2}$ but
$\arf(\uf{4})=\frac18\neq\arf(\rof{3}\otimes\uf{4})=\frac38$.

But just make any choice of $C_i$ and $R_j$ and then compute.
The only problem occurs for terms of the form $C_i\otimes R_i$
(same subscripts):
if $N_1+N_3+N_5+N_7$ for $C_i$ is not $0$ and if $N_1+N_3+N_5+N_7$
for $R_i$ is not $0$, then $B\otimes\q$ will not be tame.
Otherwise $B\otimes\q$ will be tame both (\gtpLink) and $(\gtpLink)^\prime$
hold without change.
This generalizes the Kirby\&Melvin calculation: their $\q=\uf{2}$ so
$R_1=\rof1$ and $C=C_0\perp 2C_1\perp\cdots$.
For $R_1$, $N_1+N_3+N_5+N_7=1$ so to be tame we must
have $N_1+N_3+N_5+N_7=0$ for $C_1$, which translates to no $\rof2$'s
in $C$.
Moreover $\sig_2^\prime(R_1)=1$ and
$\sig_2^\prime(C_0)=n_1-n_3+4\epsilon_H$.

\section{Topology applications.}\resetSectionLabel
\begin{ThmS}{\nextLabel[xxAd]{Theorem}}
(A theorem of Morita \cite{bMo})
Let $X$ be a $4k$ dimensional, oriented, connected
Poincar\'e duality space without boundary.
The \Pont\ square 
\[\p\colon H^{2k}(X;\cy2)\to\cy4\]
is a quadratic enhancement of the cup product pairing
$H^{2k}(X;\cy2)\times H^{2k}(X;\cy2)\to \cy2$.
Brown conjectured and Morita proved that
\[\arf(\p)=\frac{\sig(X)}8\in\qmodz\]
where $\sig(X)$ denotes the signature of $X$.
\end{ThmS}

For the proof, let $K$ denote the image of the torsion in
$H^{2k}(X;\Z)$ in $H^{2k}(X;\Z)\otimes\cy2$.
An observation of Massey \cite{bM} which he claims was 
well{--}known at the time says that, in our
notation, 
\[K^\perp/K=
\bigl(H^{2k}(X;\Z)/torsion\bigr)\otimes\cy2\ .\]

Recall that on 
$H^{2k}(X;\Z)\otimes\cy2\subset H^{2k}(X;\cy2)$
the \Pont\ square is just the cup product square reduced \MOD{4}.
This shows that
$\p$ vanishes on $K$ and that the enhancement induced by
$\p$ on $K^\perp/K$ is the \MOD{4} quadratic enhancement 
of the cup product form $B$ on $H^{2k}(X;\Z)/torsion$.
Now \xxAcxLink\ completes the proof.\qed

\begin{ThmS}{\nextLabel[xxAe]{Theorem}}
Suppose $X$ is a $4k$ dimensional, oriented, connected
Poincar\'e duality space with boundary.
The \Pont\ square 
\[\p\colon H^{2k}(X,\partial X;\cy2)\to\cy4\]
is a quadratic enhancement of the cup product pairing.
If $H^{2k}(\partial X;\Z)$ is torsion free, then
\[\arf(\p)=\frac{\sig(X)}8\in\qmodz\]
\end{ThmS}

\bigskip
Several differences arise in the bounded case.
The first is that the cup product pairing has an annihilator:
if $T=H^{2k}(X,\partial X;\cy2)$, $T^\perp$ is
the image of $H^{2k-1}(\partial X;\cy2)$
in $H^{2k}(X,\partial X;\cy2)$.
A theorem of Thomas \cite{bET} calculates that the compositions

{\small \noindent$\begin{matrix}%
H^{2k-1}(\partial X;\cy2)&\to H^{2k}(X,\partial X;\cy2)\
\xrightarrow{\p}\  H^{4k}(X,\partial X;\cy4)=&\cy4\hfill\cr\noalign{\vskip4pt}
H^{2k-1}(\partial X;\cy2)&\ \xrightarrow{\hskip 20pt x\cupp Sq^1x\hskip 10pt }\
H^{4k-1}(\partial X;\cy2)=\cy2\subset&\cy4\hfill\cr
\end{matrix}$}

\noindent are equal.
Note that if $H^{2k}(\partial X;\Z)$ is torsion{--}free, $Sq^1x=0$ so $\p$
is tame in this case.
(But not in general: any oriented boundary for an $RP^{2k-1}$
is going to have $\p$ not tame.)
We can identify $T/T^\perp$ with the image of
$H^{2k}(X,\partial X;\cy2)$ in $H^{2k}(X;\cy2)$.

Let $I_X$ denote the image of $H^{2k}(X,\partial X;\Z)$
in $H^{2k}(X;\Z)$.
Since $H^{2k}(\partial X;\Z)$ is torsion{--}free,
$I_X\otimes\cy2\subset T/T^\perp$.
Let $K\subset T/T^\perp$ denote the image of the torsion
of $I_X$ in $T/T^\perp$.

A straightforward generalization of Massey's observation
shows that $K^\perp=I_X$.
Let $B_X$ denote the matrix for the form on $I_X/torsion$
so \ $\sig(X)=\sig(B_X)$.
Since $H^{2k}(\partial X;\Z)$ is torsion{--}free, $\det B_X=\pm1$
and just as in \xxAdLink\ $\arf(\p)=\arf(\q_{B_X,1})$.
\xxAcxLink\ completes the proof.\qed

\bigskip
The function $x\cup Sq^1x$ arises in other contexts.
It is the squaring homomorphism associated to the bilinear form
$$R\colon H^k(M;\cy2)\times H^k(M,\cy2)\to\cy2$$
defined by $R(x,y)=\rof{Sq^1x\cup y,[M]}$ on any oriented
$2k+1$ dimensional manifold.
The form $R$ is symmetric.

Let $M^3$ be a closed $3$ manifold with a fixed Spin structure.
In \cite{bKT}*{p.209}, we showed how to quadratically 
enhance the linking form $\ell$ on $H^2(M;\Z)$.
Let $\q$ denote this enhancement.
We further showed that
$\arf(\q)$ is just Rochlin's $\mu${--}invariant of $M$ \MOD{8}.
Spin structures on $M$ are acted on by $H^1(M;\cy2)$.
For $x\in H^1(M;\cy2)$, define $\hat\mu(x)$ to be the difference of
the $\mu$ invariant for $M$ with its given Spin structure minus the
$\mu$ invariant for $M$ with Spin structure obtained by acting via $x$,
all reduced \MOD{8}.
We showed that $\hat\mu$ is a quadratic enhancement of $R$.
(See the remark after formula 4.9, p.213 \cite{bKT}).

Now in general $R$ is singular: in fact $R^\perp$ is precisely the 
kernel of $Sq^1$. 
We see that the image of $H^1(M;\Z)$ in $H^1(M;\cy2)$ acts trivially
on the quadratic enhancement on $\ell$ so we get
an action of $H^1(M;\cy2)/H^1(M;\Z)$ on this set.
This group is identified via the integral Bockstein $\delta$ with
$\preT 2H^2(M;\Z)$ and $R$ naturally induces a bilinear form
on $\preT 2H^2(M;\Z)$:
the enhancement $\hat\mu$ also extends to $\preT 2H^2(M;\Z)$.

In (\dmdLink) we remarked that $\preT2{H^2(M;\Z)}$ 
acts on the quadratic
enhancements of the linking form, and a comparison of (\dmdLink)
and the enhancement $\hat\mu$ shows that
$\hat\mu(x)=\mymu(\delta x)$ for any $x\in H^1(M;\cy2)$.
So the quadratic enhancement on $R$
is just the quadratic enhancement of the linking form restricted to
$\preT2{H^2(M;\Z)}$.

\makeatletter\makeatother
\begin{ThmS}{\nextLabel[xxAf]{Theorem}}
If the torsion subgroup of $H_1(M^3;\Z)$ 
is a $\cy2$ vector space, then the $\mu$ invariant \MOD{8}
is $\arf(\hat\mu)$ since $\hat\mu=\q$.
\end{ThmS}

\begin{DefS}{\nextLabel[xxAff]{Remark}}
In general, $\hat\mu$ is tame unless $\preT4(H^2(M^3;\Z))$
contains an $x$ with $\ell(x,x)=\pm\frac14$,
in which case $\hat\mu$ is not tame: e.g. the lens spaces $L(4,\pm1)$.
\end{DefS}
\end{document}